\documentclass[12pt]{amsart}
\usepackage{amsfonts}
\usepackage{amsmath,amssymb,amsthm}

\newcommand{\R}{\mathbb R}

\newcommand{\E}{\mathbb E}

\renewcommand{\span}{\mathrm{span}}
\newcommand{\tr}{\mathrm{tr}}

\newtheorem{thm}{Theorem}[section]

\newtheorem{prop}[thm]{Proposition}
\theoremstyle{definition}

\theoremstyle{remark}

\newcommand{\ds}{\displaystyle}

\begin{document}

\title[Chen Rotational Surfaces of Hyperbolic or Elliptic  Type in $\R^4_1$]
{Chen Rotational Surfaces of Hyperbolic or Elliptic  Type  in the Four-dimensional
Minkowski Space}

\author{Georgi Ganchev and Velichka Milousheva}
\address{Bulgarian Academy of Sciences, Institute of Mathematics and Informatics,
Acad. G. Bonchev Str. bl. 8, 1113 Sofia, Bulgaria}
\email{ganchev@math.bas.bg}
\address{Bulgarian Academy of Sciences, Institute of Mathematics and Informatics,
Acad. G. Bonchev Str. bl. 8, 1113, Sofia, Bulgaria; and ''L. Karavelov''
Civil Engineering Higher School, 175 Suhodolska Str., 1373 Sofia,
Bulgaria} \email{vmil@math.bas.bg}



\begin{abstract}
We study the class of spacelike surfaces in the four-dimensional
Minkowski space whose mean curvature vector at any point is a
non-zero spacelike vector or  timelike vector. These surfaces are
determined up to a motion by eight invariant functions satisfying
some natural conditions. The subclass of Chen surfaces is
characterized by the condition  one of these invariants to be
zero. In the present paper we describe all Chen spacelike
rotational surfaces of hyperbolic or elliptic type.

\end{abstract}

\maketitle

\section{Preliminaries}

In \cite{GM5} we considered the general theory of spacelike surfaces
in the four-dimensional Minkowski space $\R^4_1$. The basic
feature of our treatment of these surfaces was the introduction of an
invariant linear map of Weingarten-type in the tangent plane at
any point of the surface, following the  approach to the theory of surfaces in $\R^4$ \cite{GM1,GM4}.
 Studying surfaces in the Euclidean space
$\R^4$, in \cite{GM1} we introduced a linear map $\gamma$ of
Weingarten-type, which plays a similar role in the theory of
surfaces in $\R^4$ as the Weingarten map in the theory of surfaces
in $\R^3$. The map $\gamma$ generates the corresponding second
fundamental form $II$ at any point $p \in M^2$ in the standard
way. We gave a geometric interpretation of the second
fundamental form and the Weingarten map of the surface in
\cite{GM3}.

 Let
$M^2: z=z(u,v), \,\, (u,v) \in \mathcal{D}$ $(\mathcal{D} \subset \R^2)$  be a spacelike surface in $\R^4_1$ with
tangent space $T_pM^2 = \span \{z_u,z_v\}$ at an arbitrary point $p=z(u,v)$ of $M^2$. Since $M^2$ is spacelike,
$\langle z_u,z_u \rangle > 0$, $\langle z_v,z_v \rangle > 0$.
We use the standard denotations
$E(u,v)=\langle z_u,z_u \rangle, \; F(u,v)=\langle z_u,z_v
\rangle, \; G(u,v)=\langle z_v,z_v \rangle$ for the coefficients
of the first fundamental form
$I(\lambda,\mu):= E \lambda^2 + 2F \lambda \mu + G \mu^2,\,\,
\lambda, \mu \in \R$. Since $I(\lambda, \mu)$ is positive
definite we set $W=\sqrt{EG-F^2}$.
We choose a normal frame field $\{n_1, n_2\}$ such that $\langle
n_1, n_1 \rangle =1$, $\langle n_2, n_2 \rangle = -1$, and the
quadruple $\{z_u,z_v, n_1, n_2\}$ is positively oriented in
$\R^4_1$.

Considering the tangent
space $T_pM^2$ at a point $p \in M^2$, in \cite{GM5} we introduced an invariant
$\zeta_{\,g_1,g_2}$ of a pair of two tangents $g_1$, $g_2$ using
the second fundamental tensor $\sigma$ of $M^2$. By means of this
invariant we defined conjugate, asymptotic, and principal
tangents.
The second fundamental form $II$ of the surface $M^2$ at a point
$p \in M^2$ is introduced on the base of conjugacy of two tangents
at the point.
The coefficients $L, M, N$ of the second fundamental form $II$ are determined as follows:
$$L = \ds{\frac{2}{W}} \left|%
\begin{array}{cc}
\vspace{2mm}
  c_{11}^1 & c_{12}^1 \\
  c_{11}^2 & c_{12}^2 \\
\end{array}%
\right|; \quad
M = \ds{\frac{1}{W}} \left|%
\begin{array}{cc}
\vspace{2mm}
  c_{11}^1 & c_{22}^1 \\
  c_{11}^2 & c_{22}^2 \\
\end{array}%
\right|; \quad
N = \ds{\frac{2}{W}} \left|%
\begin{array}{cc}
\vspace{2mm}
  c_{12}^1 & c_{22}^1 \\
  c_{12}^2 & c_{22}^2 \\
\end{array}%
\right|,$$
where the functions $c_{ij}^k, \,\, i,j,k = 1,2$ are given by
$$\begin{array}{lll}
\vspace{2mm}
c_{11}^1 = \langle z_{uu}, n_1 \rangle; & \qquad   c_{12}^1 = \langle z_{uv}, n_1 \rangle; & \qquad  c_{22}^1 = \langle z_{vv}, n_1 \rangle;\\
\vspace{2mm}
c_{11}^2 = \langle z_{uu}, n_2 \rangle; & \qquad  c_{12}^2 = \langle z_{uv}, n_2 \rangle; & \qquad
c_{22}^2 = \langle z_{vv}, n_2 \rangle.
\end{array} $$

The second fundamental form $II$ determines an invariant linear
map $\gamma$  of Weingarten-type at any point of the surface,
which generates  two invariant functions:
$$k := \det \gamma = \frac{LN - M^2}{EG - F^2}, \qquad
\varkappa := -\frac{1}{2}\,{\rm tr}\, \gamma =
\frac{EN+GL-2FM}{2(EG-F^2)}.$$
The functions $k$ and $\varkappa$
are invariant under changes of the parameters of the surface and
changes of the normal frame field.
 The sign of $k$ is invariant under congruences
 and the sign of $\varkappa$ is invariant under motions
in $\R^4_1$. However, the sign of $\varkappa$ changes under
symmetries with respect to a hyperplane in $\R^4_1$. We proved that
the invariant $\varkappa$ is the curvature of the normal
connection of the surface. The number of asymptotic tangents at a
point of $M^2$ is determined by the sign of the invariant $k$. In
the case $k=0$ there exists a one-parameter family of asymptotic
lines, which are principal.

It is interesting to note that the
''umbilical'' points, i.e. points at which the coefficients of the
first and the second fundamental forms are proportional, are
exactly the points at which the mean curvature vector $H$ is zero.
So, the spacelike surfaces consisting of ''umbilical'' points in
$\R^4_1$ are exactly the minimal surfaces. Minimal spacelike
surfaces are characterized in terms of the invariants $k$ and
$\varkappa$ by the equality $\varkappa^2 - k =0$.

Analogously to $\R^3$ and $\R^4$, the invariants $k$ and
$\varkappa$ divide the points of $M^2$ into four types: flat,
elliptic, hyperbolic and parabolic points. The surfaces consisting
of flat points  are characterized by the conditions $k= \varkappa
= 0$. We gave a local geometric description of spacelike surfaces
 consisting of flat points whose mean curvature vector
at any point is a non-zero spacelike vector or  timelike vector,
proving that any such a surface either lies in a hyperplane of
$\R^4_1$ or is part of a developable ruled surface in $\R^4_1$ \cite{GM5}.

Using the introduced principal lines on a spacelike surface in
$\R^4_1$ whose mean curvature vector at any point is a non-zero
spacelike vector or  timelike vector, we found a geometrically
determined moving frame field  on such a surface. Writing the
derivative formulas of Frenet-type for this frame field, we
obtained eight invariant functions and proved a fundamental
theorem of Bonnet-type, stating that these eight invariants under
some natural conditions determine the surface up to a motion in
$\R^4_1$.

One of these eight invariants is closely related to the theory of Chen surfaces.
We shall recall the notion of Chen submanifolds.
Let $M^n$ be an $n$-dimensional submanifold of
$(n+m)$-dimensional Riemannian manifold $\widetilde{M}^{n+m}$ and
$\xi$ be a normal vector field of $M^n$. In \cite{Chen1} B.-Y.
Chen defined the \emph{allied vector field} $a(\xi)$ of $\xi$ by
the formula
$$a(\xi) = \ds{\frac{\|\xi\|}{n} \sum_{k=2}^m \{\tr(A_1 A_k)\}\xi_k},$$
where $\{\xi_1 = \ds{\frac{\xi}{\|\xi\|}},\xi_2,\dots, \xi_m \}$
is an orthonormal base of the normal space of $M^n$, and $A_i =
A_{\xi_i}, \,\, i = 1,\dots, m$ is the shape operator with respect
to $\xi_i$. In particular, the allied vector field $a(H)$ of the
mean curvature vector field $H$ is a well-defined normal vector
field which is orthogonal to $H$. It is called the \emph{allied
mean curvature vector field} of $M^n$ in $\widetilde{M}^{n+m}$.
B.-Y. Chen defined  the $\mathcal{A}$-submanifolds to be those
submanifolds of $\widetilde{M}^{n+m}$ for which
 $a(H)$ vanishes identically \cite{Chen1}.
In \cite{GVV1,GVV2} the $\mathcal{A}$-submanifolds are called
\emph{Chen submanifolds}. It is easy to see that minimal
submanifolds, pseudo-umbilical submanifolds and hypersurfaces are
Chen submanifolds. These Chen submanifolds are said to be \emph{trivial
$\mathcal{A}$-submanifolds}.

Now, let $M^2$ be a spacelike surface in  $\R^4_1$ whose mean curvature vector at any point is a non-zero
spacelike vector or  timelike vector.
We denote by $x$ and $y$  the principal unit tangent vector fields of $M^2$, and by $H$ - the mean curvature vector field.
One of the invariants in the Frenet-type derivative formulas of  $M^2$ (eg. \cite{GM5}) is
 $\lambda = \ds{\frac{1}{\sqrt{\langle H,H \rangle}}} \, \langle\sigma(x,y), H \rangle$ in the case
 $\langle H,H \rangle >0$, and
 $\lambda = \ds{- \frac{1}{\sqrt{- \langle H,H \rangle}}} \, \langle\sigma(x,y), H \rangle$ in the case
 $\langle H,H \rangle <0$.
Applying the definition of the allied mean curvature vector field and the derivative formulas of $M^2$ we get
$$a(H) = \ds{\frac{\sqrt{\varkappa^2-k}}{2} \, \lambda \,l}.$$
Hence, if $M^2$ is free of minimal points ($\varkappa^2-k \neq 0$), then $a(H) = 0$ if and
only if $\lambda = 0$. This gives the geometric interpretation of the
invariant $\lambda$. It is clear that $M^2$ is a non-trivial  Chen
surface if and only if the invariant $\lambda$ is zero.

\vskip 2mm
In the present paper we study spacelike rotational  surfaces of hyperbolic or elliptic type
in the four-dimensional Minkowski space $\R^4_1$ and
we describe the class of Chen  rotational surfaces of hyperbolic or elliptic type.

\section{Rotational surfaces with two-dimensional axis in $\R^4_1$}

In \cite{GM1} we considered the class of the rotational surfaces with two-dimensional axis
in the four-dimensional Euclidean space $\R^4$.

Let $Oe_1e_2e_3e_4$ be a fixed orthonormal base of $\R^4$ and
$\R^3$ be the subspace spanned by $e_1, e_2, e_3$. We consider a
smooth curve $c: \widetilde{z} = \widetilde{z}(u), \,\, u \in J$
in $\R^3$, parameterized by
$$\widetilde{z}(u) = \left( x_1(u), x_2(u), r(u)\right); \quad u \in J.$$
Without loss of generality we assume that $c$ is
parameterized by the arc-length, i.e. $(x_1')^2 + (x_2')^2 +
(r')^2 = 1$. We assume also that $r(u)>0,\,\, u \in J$.
Let $\kappa$ and $\tau$ be the curvature and the torsion of  $c$.
We denote by $c_1$ the projection of $c$ into the 2-dimensional
plane $Oe_1e_2$ and by  $\kappa_1$ the curvature of the plane curve $c_1$, i.e.
$\kappa_1 = x_1' x_2'' - x_1'' x_2'$.

Let us consider the rotational surface $M^2$ in $\R^4$ given by
$$z(u,v) = \left( x_1(u), x_2(u), r(u) \cos v, r(u) \sin v\right);
\quad u \in J,\,\,  v \in [0; 2\pi).$$
$M^2$ is obtained by the rotation of the curve $c$ about the two-dimensional axis $Oe_1e_2$
(the rotation of $c$ that leaves the plane $Oe_1e_2$ fixed).

In \cite{GM1} we found that the invariants $k$, $\varkappa$ and the Gauss curvature $K$ of the
rotational surface $M^2$ are expressed as follows:
$$k = - \displaystyle{\frac{(\kappa_1)^2}{r^2}}; \quad \quad \varkappa = 0; \quad \quad
K = - \displaystyle{\frac{r''}{r}}.$$

Obviously, the rotational surface $M^2$ is a surface with flat normal connection since $\varkappa = 0$.
We described all
rotational surfaces, for which the invariant $k$ is constant \cite{GM1}.

The mean curvature vector field $H$ of $M^2$ is given by:
$$H = \ds{\frac{\kappa^2 r - r''}{2\kappa r}\, e_1 - \frac{\kappa_1}{2 \kappa r}\, e_2,}$$
where $e_1$ and $e_2$ are the normal vector fields, defined by
$$\begin{array}{l}
\vspace{2mm}
e_1 = \displaystyle{\frac{1}{\kappa}\left(x_1'', x_2'', r'' \cos v, r'' \sin v \right)};\\
\vspace{2mm} e_2 = \displaystyle{\frac{1}{\kappa} \left( x_2' r''
- x_2'' r' , x_1'' r' - x_1' r'', (x_1' x_2'' - x_1'' x_2') \cos v, (x_1' x_2'' - x_1''
x_2') \sin v \right)}.
\end{array}$$

Hence, the mean curvature vector field $H$ vanishes if and only if $\kappa_1 =0$ and $\kappa^2 r - r'' = 0$.
In this case the rotational surface $M^2$ is a trivial Chen surface ($M^2$ is minimal).

In the case $\kappa_1 =0$ and $\kappa^2 r - r'' \neq 0$ the rotational surface
$M^2$ lies in a three-dimensional subspace $\R^3$ of $\R^4$.
Moreover, in this case one can easily get that the allied vector field $a(H)$ of the mean
curvature vector field  is zero, and hence $M^2$ is a  Chen  surface in $\R^3$. Hence, $M^2$ is a trivial Chen surface.

In the case when $\kappa_1 \neq 0$ the rotational surface $M^2$ is a non-minimal surface in $\R^4$ free of flat points.
Calculating the invariant $\lambda$ of $M^2$ we get that  $\lambda = 0$ if and only if
$$\kappa^4 r^2 - (r'')^2 - (\kappa_1)^2 = 0.$$
In this case $M^2$ is a non-trivial Chen surface in $\R^4$.

\vskip 3mm
Now we shall study spacelike rotational  surfaces of hyperbolic or elliptic type
in Minkowski space $\R^4_1$ and we shall describe the class of Chen  rotational surfaces.

We consider the Minkowski space $\R^4_1$ endowed with the metric
$\langle , \rangle$ of signature $(3,1)$.
 Let $Oe_1e_2e_3e_4$ be a
fixed orthonormal coordinate system in $\R^4_1$, i.e. $e_1^2 =
e_2^2 = e_3^2 = 1, \, e_4^2 = -1$, giving the orientation of
$\R^4_1$. The standard flat metric is given in local coordinates by
$dx_1^2 + dx_2^2 + dx_3^2 -dx_4^2.$

A surface $M^2: z = z(u,v), \, \, (u,v) \in {\mathcal D}$
(${\mathcal D} \subset \R^2$) in $\R^4_1$ is said to be
\emph{spacelike} if $\langle , \rangle$ induces  a Riemannian
metric $g$ on $M^2$. Thus at each point $p$ of a spacelike surface
$M^2$ we have the following decomposition $\R^4_1 = T_pM^2 \oplus
N_pM^2$ with the property that the restriction of the metric
$\langle , \rangle$ onto the tangent space $T_pM^2$ is of
signature $(2,0)$, and the restriction of the metric $\langle ,
\rangle$ onto the normal space $N_pM^2$ is of signature $(1,1)$.

We  consider a
smooth spacelike  curve $c: \widetilde{z} = \widetilde{z}(u), \,\, u \in J$, parameterized by
$$\widetilde{z}(u) = \left( x_1(u), x_2(u), 0, r(u)\right); \quad u \in J.$$
The curve $c$ lies in the three-dimensional subspace $\R^3_1 = \span\{e_1, e_2, e_4\}$ of $\R^4_1$.
Without loss of generality we assume that $c$ is
parameterized by the arc-length, i.e. $(x_1')^2 + (x_2')^2 - (r')^2 = 1$.
We assume also that  $r(u)>0, \,\, u \in J$ and $\langle t_c'(u), t_c'(u) \rangle \neq 0,\,\,u \in J$,
where $t_c(u)=z'(u)$.
We have the following possibilities:

1.) $\langle t_c', t_c'\rangle > 0$, i.e. $(x_1'')^2 + (x_2'')^2 - (r'')^2 > 0$.
In this case the Frenet formulas of $c$ are given by
$$t'_c = \kappa\,n_c; \quad n'_c = - \kappa\,t_c + \tau\,b_c; \quad b'_c = - \tau\, n_c,$$
where $\{t_c, n_c, b_c\}$ is the Frenet frame field of $c$; $\kappa$ and $\tau$ - the curvature and the torsion of $c$;
$n_c$ is  spacelike, $b_c$ is timelike.

\vskip 1mm
2.) $\langle t_c', t_c' \rangle < 0$, i.e. $(x_1'')^2 + (x_2'')^2 - (r'')^2 < 0$.
In this case we have the following Frenet formulas of $c$:
$$t'_c = \kappa\,n_c; \quad n'_c = \kappa\,t_c + \tau\,b_c; \quad b'_c = \tau\, n_c,$$
where $n_c$ is  timelike, $b_c$ is spacelike.

\vskip 2mm
Let us consider the surface $M^2$ in $\R^4_1$ given by
$$z(u,v) = \left( x_1(u), x_2(u), r(u) \sinh v, r(u) \cosh v\right);
\quad u \in J,\,\,  v \in \R.\leqno{(2.1)}$$
The tangent space of $M^2$ is spanned by the vector fields
$$\begin{array}{l}
\vspace{2mm}
z_u = \left(x_1', x_2', r' \sinh v, r' \cosh v \right);\\
\vspace{2mm} z_v = \left(0, 0, r \cosh v, r \sinh v \right).
\end{array}$$
Hence, the coefficients of the first fundamental form of $M^2$ are:
$$E = \langle z_u, z_u \rangle = 1; \quad F = \langle z_u, z_v \rangle = 0; \quad G = \langle z_v, z_v \rangle = r^2(u),$$
and the induced metric $g$ on $M^2$ is a Riemannian metric:
$$g = du^2 + r^2(u) dv^2.$$
So, the surface $M^2$ in $\R^4_1$, defined by (2.1), is a spacelike surface.
It is called  a \emph{spacelike rotational surface of hyperbolic type} in $\R^4_1$.
It is an orbit of a spacelike regular curve under the action of the orthogonal transformations of $\R^4_1$
which leave a spacelike plane point-wise fixed. In our case the two-dimensional plane $Oe_1e_2$ is fixed.

A classification of all timelike and spacelike hyperbolic rotational surfaces
with non-zero constant mean curvature in the three-dimensional de Sitter space
$\mathbb{S}^3_1$ is given in \cite{Liu-Liu}.
Similarly, a classification of the spacelike and timelike Weingarten rotation surfaces in $\mathbb{S}^3_1$
is found in \cite{Liu-Liu-2}.

Here we shall describe the class of  Chen spacelike rotational surfaces of hyperbolic type in $\R^4_1$.

The second
partial derivatives of $z(u,v)$ are expressed as follows
$$\begin{array}{l}
\vspace{2mm}
z_{uu} = \left(x_1'', x_2'', r'' \sinh v, r'' \cosh v \right);\\
\vspace{2mm}
z_{uv} = \left(0, 0, r' \cosh v, r' \sinh v \right);\\
\vspace{2mm} z_{vv} = \left(0, 0, r \sinh v, r \cosh v \right).
\end{array}$$
We consider the following orthonormal tangent vector fields
$$\begin{array}{l}
\vspace{2mm}
\overline{x} = \left(x_1', x_2', r' \sinh v, r' \cosh v \right);\\
\vspace{2mm} \overline{y} = \left(0, 0,  \cosh v, \sinh v \right),
\end{array}$$
i.e. $z_u = \overline{x}, \,\, z_v = r\, \overline{y}$, and the
normal vector fields $n_1$, $n_2$, defined by
$$\begin{array}{l}
\vspace{2mm}
n_1 = \displaystyle{\frac{1}{\kappa}\left(x_1'', x_2'', r'' \sinh v, r'' \cosh v \right)};\\
\vspace{2mm}
n_2 = \displaystyle{\frac{1}{\kappa}
\left( x_2' r'' - x_2'' r' , r' x_1''  - x_1' r'', (x_2' x_1'' - x_1' x_2'') \sinh
v, (x_2' x_1'' - x_1' x_2'') \cosh v \right)}.
\end{array}$$

In the case when $n_c$ is spacelike the normal vector field $n_1$ is spacelike and $n_2$ is timelike.
If $n_c$ is timelike, then  $n_1$ is timelike and $n_2$ is spacelike.
We denote $\varepsilon = \langle n_c, n_c \rangle$. Then we have $\langle n_1, n_1 \rangle = \varepsilon$,
$\langle n_2, n_2 \rangle = -\varepsilon$.

We calculate the coefficients of the second fundamental tensor of $M^2$:
$$\begin{array}{ll}
\vspace{2mm}
c_{11}^1 = \langle z_{uu}, n_1 \rangle = \varepsilon \,\kappa; & \qquad c_{11}^2 = \langle z_{uu}, n_2 \rangle = 0;\\
\vspace{2mm}
c_{12}^1 = \langle z_{uv}, n_1 \rangle = 0; & \qquad c_{12}^2 = \langle z_{uv}, n_2 \rangle = 0;\\
\vspace{2mm} c_{22}^1 = \langle z_{vv}, n_1 \rangle = \ds{-\frac{r
r''}{\kappa}}; & \qquad c_{22}^2 = \langle z_{vv}, n_2 \rangle =
\ds{\frac{r}{\kappa} (x_1' x_2'' - x_2' x_1'')}.
\end{array} \leqno{(2.2)}$$

Let us denote by $c_1$ the projection of $c$ into the 2-dimensional
plane $Oe_1e_2$ and by  $\kappa_1$ the curvature of $c_1$, i.e.
$\kappa_1 = x_1' x_2'' - x_2' x_1''$.

Using (2.2) we calculate the coefficients of the second fundamental form of $M^2$:
$$L = 0; \quad M = \varepsilon \kappa_1; \quad N = 0.$$
Hence, the invariants $k$ and $\varkappa$ of the rotational surface of hyperbolic  type are expressed  as
$$k = - \displaystyle{\frac{(\kappa_1)^2}{r^2}}; \quad \quad \varkappa = 0.$$

Consequently, any spacelike  rotational surface of hyperbolic  type in $\R^4_1$ is a
surface with flat normal connection since $\varkappa = 0$.

With respect to the frame field
$\{\overline{x}, \overline{y}, n_1, n_2\}$ the derivative formulas
of $M^2$ look like:
$$\begin{array}{ll}
\vspace{2mm}
\nabla'_{\overline{x}}\overline{x} = \quad \quad \quad \quad \quad  \kappa\,n_1 ;  & \qquad
\nabla'_{\overline{x}}n_1 = - \varepsilon \kappa\,\overline{x}\quad
\quad \quad \quad \quad \quad +\tau\,n_2;\\
\vspace{2mm}
\nabla'_{\overline{x}}\overline{y} = 0;  & \qquad
\nabla'_{\overline{y}} n_1 = \quad  \quad \quad
\displaystyle{\frac{r''}{r \kappa}}\,\overline{y} ;\\
\vspace{2mm}
\nabla'_{\overline{y}}\overline{x} = \quad
\quad\;\displaystyle{\frac{r'}{r}}\,\overline{y} ;  & \qquad
\nabla'_{\overline{x}}n_2 = \quad \quad \quad \quad \quad\quad -\varepsilon \tau\,n_1;\\
\vspace{2mm}
\nabla'_{\overline{y}}\overline{y} =
-\displaystyle{\frac{r'}{r}}\,\overline{x} \quad \;-
\displaystyle{\frac{r''}{r \kappa}}\,\varepsilon n_1 -
\displaystyle{\frac{\kappa_1}{r \kappa}}\,\varepsilon n_2; & \qquad
\nabla'_{\overline{y}}n_2 = \quad \quad  \quad
\displaystyle{- \frac{\kappa_1}{r \kappa}}\,\overline{y}.
\end{array} \leqno{(2.3)}$$
So, the Gauss curvature of $M^2$ is:
$$K = - \displaystyle{\frac{r''}{r}}.$$

Obviously $M^2$ is not parameterized by the principal
lines. The principal tangents of $M^2$ are:
$$x = \displaystyle{\frac{\sqrt{2}}{2}\, \overline{x} + \frac{\sqrt{2}}{2} \, \overline{y}}; \qquad
y = \displaystyle{\frac{\sqrt{2}}{2}\, \overline{x} -
\frac{\sqrt{2}}{2} \, \overline{y}}. \leqno{(2.4)}$$
Then formulas (2.3) and (2.4) imply that the normal vector fields $\sigma(x,x)$, $\sigma(x,y)$, and $\sigma(y,y)$ are given by
$$\begin{array}{l}
\vspace{2mm}
\sigma(x,x)=\ds{\frac{r \kappa^2 - \varepsilon r''}{2r \kappa} \, n_1 - \frac{\varepsilon  \kappa_1}{2r \kappa}\, n_2},\\
\vspace{2mm}
\sigma(x,y)= \ds{\frac{r \kappa^2 + \varepsilon r''}{2r \kappa} \, n_1 +\frac{\varepsilon  \kappa_1}{2r \kappa}\,  n_2},\\
\vspace{2mm}
\sigma(y,y) =\ds{\frac{r \kappa^2 - \varepsilon r''}{2r \kappa} \, n_1 - \frac{\varepsilon  \kappa_1}{2r \kappa}\, n_2}.
\end{array} \leqno{(2.5)}$$
The normal  mean curvature vector field of $M^2$  is
$$ H = \ds{\frac{r \kappa^2 - \varepsilon r''}{2r \kappa} \, n_1 - \frac{\varepsilon  \kappa_1}{2r \kappa}\, n_2}. \leqno{(2.6)}$$

All Chen spacelike rotational surfaces of hyperbolic  type in $\R^4_1$ are described in the following

\begin{prop} \label{P:hyperbolic rotation}
The spacelike rotational surface  of hyperbolic  type $M^2$, defined by (2.1), is
 a Chen surface if and only if one of the following cases holds:

(i) $\kappa_1 =0$ and $r \kappa^2 - \varepsilon r'' = 0$; in such case $M^2$ is a minimal surface (trivial Chen surface);

(ii) $\kappa_1 =0$ and $r \kappa^2 - \varepsilon  r'' \neq 0$; in such case
$M^2$ is a surface lying in a three-dimensional subspace of $\R^4_1$ (trivial Chen surface);

(iii) $\kappa_1 \neq 0$ and $r^2 \kappa^4  - (r'')^2 + (\kappa_1)^2 = 0$; in such case
$M^2$ is a non-trivial Chen surface in $\R^4_1$.
\end{prop}

\proof
\emph{(i)} From (2.6) it follows that the mean curvature vector field $H$ vanishes if and only if
$$\kappa_1 =0; \qquad r \kappa^2  - \varepsilon  r'' = 0.$$
In this case the rotational surface $M^2$ is a trivial Chen surface ($M^2$ is minimal).

\emph{(ii)} In the case $\kappa_1 =0$ and $r \kappa^2  - \varepsilon  r'' \neq 0$ we get that
$k = \varkappa =0$, and hence $M^2$ consists of flat points. We shall prove that $M^2$ lies in a hyperplane of
$\R^4_1$. The equality $\kappa_1 =0$ implies that the projection of the curve $c$ into  $Oe_1e_2$ lies on a straight line, and
hence $c$ lies in a two-dimensional plane orthogonal to $Oe_1e_2$. Thus the torsion of $c$ is $\tau = 0$.
Then, using derivative formulas (2.3) we get
$$\nabla'_{\overline{x}}n_2 = 0; \qquad \nabla'_{\overline{y}}n_2 = 0,$$
which imply that the normal vector field $n_2$ is constant.
Consequently, $M^2$ lies in the hyperplane $\E^3$ of $\R^4_1$ orthogonal to $n_2$, i.e. $\E^3 = \span \{x, y, n_1\}$.
Hence, $M^2$ is a trivial Chen surface.

 \emph{(iii)}
In the case when $\kappa_1 \neq 0$ the rotational surface $M^2$ is a non-minimal surface in $\R^4$ free of flat points.
According to \cite{GM5} $M^2$ is a Chen surface if and only if the invariant $\lambda$ is zero.
We have that  $\lambda = 0$ if and only if
$\langle \sigma(x,y), H \rangle = 0$,  which in view of (2.5) and (2.6) implies
$$r^2 \kappa^4  - (r'')^2 + (\kappa_1)^2 = 0.$$
In this case $M^2$ is a non-trivial Chen surface in $\R^4_1$.

\qed

\vskip 3mm
In a similar way we consider a spacelike surface in $\R^4_1$ which is an orbit of
a spacelike regular curve $c$ under the action of the orthogonal transformations of $\R^4_1$
which leave a timelike plane point-wise fixed.
Let us consider a spacelike curve
$c: \widetilde{z} = \widetilde{z}(u), \,\, u \in J$, parameterized by
$$\widetilde{z}(u) = \left( r(u), 0, x_1(u), x_2(u)\right); \quad u \in J.$$
The curve $c$ lies in the three-dimensional subspace $\R^3_1 = \span\{e_1, e_3, e_4\}$ of $\R^4_1$.
Again we assume that $c$ is
parameterized by the arc-length, i.e. $(r')^2 + (x_1')^2 - (x_2')^2 = 1$.
We assume also that  $r(u)>0, \,\, u \in J$ and $\langle t_c'(u), t_c'(u) \rangle \neq 0,\,\,u \in J$,
where $t_c(u)=z'(u)$.

Now we consider the surface $M^2$ in $\R^4_1$ given by
$$z(u,v) = \left(r(u) \cos v, r(u) \sin v, x_1(u), x_2(u)\right);
\quad u \in J,\,\,  v \in [0; 2\pi).\leqno{(2.7)}$$
The tangent space of $M^2$ is spanned by the vector fields
$$\begin{array}{l}
\vspace{2mm}
z_u = \left(r' \cos v, r' \sin v, x_1', x_2' \right);\\
\vspace{2mm} z_v = \left( - r \sin v, r \cos v, 0, 0 \right).
\end{array}$$
Hence, the coefficients of the first fundamental form of $M^2$ are
$$E = \langle z_u, z_u \rangle = 1; \quad F = \langle z_u, z_v \rangle = 0; \quad G = \langle z_v, z_v \rangle = r^2(u),$$
and the induced metric $g$ on $M^2$ is a Riemannian metric:
$$g = du^2 + r^2(u) dv^2.$$
The surface $M^2$, defined by (2.7), is a spacelike surface  in $\R^4_1$.
It is obtained by the rotation of the curve $c$ about the two-dimensional Lorentz plane $Oe_3e_4$.
It is called a \emph{spacelike rotational surface of elliptic type}.
A local classification of spacelike surfaces in $\R^4_1$, which are invariant under spacelike rotations,
and with  mean curvature vector either vanishing or lightlike, is obtained in \cite{Hae-Ort}.
Here we shall describe the class of   Chen spacelike rotational surfaces of elliptic type in $\R^4_1$.

As in the case of rotational surfaces of hyperbolic  type we consider the orthonormal tangent vector fields
$$\overline{x} = z_u; \qquad \overline{y} = \ds{\frac{z_v}{r}},$$
and the normal vector fields $n_1$, $n_2$, defined by
$$\begin{array}{l}
\vspace{2mm}
n_1 = \displaystyle{\frac{1}{\kappa}\left(r'' \cos v, r'' \sin v, x_1'', x_2'' \right)};\\
\vspace{2mm}
n_2 = \displaystyle{\frac{1}{\kappa}
\left((x_1' x_2'' - x_2' x_1'') \cos v, (x_1' x_2'' - x_2' x_1'') \sin v,   x_2' r'' -  r' x_2'', x_1' r'' - r' x_1'' \right)},
\end{array}$$
where $\kappa$ is the curvature  of $c$.
The principal tangents of $M^2$ are
$x = \displaystyle{\frac{\sqrt{2}}{2} (\overline{x} + \overline{y})}; \,\,
y = \displaystyle{\frac{\sqrt{2}}{2} (\overline{x} -  \overline{y})}$.
With respect to the principal tangents $x,y$ we get the following formulas:
$$\begin{array}{l}
\vspace{2mm}
\sigma(x,x)=\ds{\frac{r \kappa^2 - \varepsilon r''}{2r \kappa} \, n_1 + \frac{\varepsilon \kappa_1}{2r \kappa}\, n_2},\\
\vspace{2mm}
\sigma(x,y)= \ds{\frac{r \kappa^2 + \varepsilon r''}{2r \kappa} \, n_1 - \frac{\varepsilon \kappa_1}{2r \kappa}\, n_2},\\
\vspace{2mm}
\sigma(y,y) =\ds{\frac{r \kappa^2 - \varepsilon r''}{2r \kappa} \, n_1 + \frac{\varepsilon \kappa_1}{2r \kappa}\, n_2}.
\end{array}$$
The invariants $k$, $\varkappa$ and the Gauss curvature $K$ of $M^2$ are expressed  as in the hyperbolic case:
$$k = - \displaystyle{\frac{(\kappa_1)^2}{r^2}}; \qquad \varkappa = 0; \qquad K = \ds{-\frac{r''}{r}}.$$
The normal  mean curvature vector field of $M^2$  is
$$ H = \ds{\frac{r \kappa^2 - \varepsilon r''}{2r \kappa} \, n_1 + \frac{\varepsilon \kappa_1}{2r \kappa}\,  n_2}.$$

All Chen spacelike rotational surfaces of elliptic type in $\R^4_1$ are described in the following

\begin{prop} \label{P:spacelike rotation}
The spacelike  rotational surface of elliptic type $M^2$, defined by (2.7), is
 a Chen surface if and only if one of the following cases holds:

(i) $\kappa_1 =0$ and $r \kappa^2 - \varepsilon r'' = 0$; in such case $M^2$ is a minimal surface (trivial Chen surface);

(ii) $\kappa_1 =0$ and $r \kappa^2 - \varepsilon  r'' \neq 0$; in such case
$M^2$ is a surface lying in a three-dimensional subspace of $\R^4_1$ (trivial Chen surface);

(iii) $\kappa_1 \neq 0$ and $r^2 \kappa^4  - (r'')^2 + (\kappa_1)^2 = 0$; in such case
$M^2$ is a non-trivial Chen surface in $\R^4_1$.
\end{prop}

The proof is similar to the proof of Proposition \ref{P:hyperbolic rotation}.

\vskip 3mm At the end of the section we shall describe all spacelike
rotational surfaces of hyperbolic or elliptic type in $\R^4_1$, for which the invariant $k$ is constant.

\vskip 1mm 1. The invariant $k = 0$ if and only if $\kappa_1 = 0$,
i.e. the projection of the curve $c$ into the two-dimensional axis lies on a straight line.
In this case $M^2$ lies in a hyperplane $\E^3$ of $\R^4_1$. There are two subcases:

\hskip 10mm 1.1. If $K = 0$, i.e. $r'' = 0$,  then $M^2$ is a
developable ruled surface in $\E^3$.

\hskip 10mm 1.2. If  $K \neq 0$, i.e. $r'' \neq 0$, then $M^2$ is
a non-flat surface in $\E^3$.

\vskip 1mm 2. The invariant $k = {\rm const}$ ($k \neq 0$) if and
only if $r(u) = a\, (x_1' x_2'' -  x_2' x_1'')$,  $a = {\rm
const}$, $a\neq 0$. Moreover, if $r(u)$ satisfies $r''(u) = c\, r(u)$, then
the Gauss curvature $K$ is also a constant.

\vskip 5mm \textbf{Acknowledgements:} The second author is
partially supported by "L. Karavelov" Civil Engineering Higher
School, Sofia, Bulgaria under Contract No 10/2010.


\begin{thebibliography}{99}

\bibitem{Chen1}
Chen B.-Y.,  Geometry of submanifolds, Marcel Dekker, Inc., New
York, 1973.

\bibitem{GM1}
Ganchev G., Milousheva V.,  On the theory of surfaces in the
four-dimensional Euclidean space,   Kodai Math. J., 2008,  \textbf{ 31},
183--198.

\bibitem{GM3}
Ganchev G., Milousheva V., Invariants of lines on surfaces in
$\R^4$,  C. R. Acad. Bulg. Sci., 2010, \textbf{63}, (6), 835--842.

\bibitem{GM4}
Ganchev G., Milousheva V., Invariants and Bonnet-type theorem for
 surfaces in $\R^4$, Cent. Eur. J. Math., 2010, \textbf{8} (6), 993--1008.

\bibitem{GM5}
Ganchev G., Milousheva V., An invariant theory of spacelike
surfaces in the four-dimensional Minkowski space,  Mediterr. J.
Math., 2010, DOI: 10.1007/s00009-010-0108-2 .

\bibitem{GVV1}
Gheysens L.,  Verheyen P.,  Verstraelen L., Sur les surfaces
$\mathcal{A}$ ou les surfaces de Chen,   C. R. Acad. Sci. Paris,
S\'{e}r. I, 1981, \textbf{292}, 913--916.

\bibitem{GVV2}
Gheysens L.,  Verheyen P.,  Verstraelen L., Characterization and
examples of Chen submanifolds,  J. Geom., 1983,  \textbf{20}, 47--62.

\bibitem{Hae-Ort}
Haesen S.,  Ortega M., Marginally trapped surfaces in Minkowski 4-space
invariant under a rotation subgroup of the Lorentz group,  Gen. Relativ. Gravit., 2009,  \textbf{41}, 1819--1834.

\bibitem{Liu-Liu}
Liu H., Liu G., Hyperbolic rotation surfaces of constant mean curvature in 3-de Sitter space,
Bull. Belg. Math. Soc., 2000, \textbf{7},
455--466.

\bibitem{Liu-Liu-2}
Liu H., Liu G., Weingarten rotation surfaces  in 3-imensional de Sitter space,
J. Geom., 2004, \textbf{79},
156--168.

\end{thebibliography}
\end{document}